\documentclass[preprint,prb,nobibnotes,showkeys,preprintnumbers,amsmath,amssymb]{revtex4}
\usepackage{graphicx}
\usepackage{subfig}
\usepackage{caption}
\usepackage{pst-all}
\usepackage{dcolumn}
\usepackage{bm}
\usepackage[dvips]{epsfig}
\usepackage{rotating}
\usepackage{bm}
\usepackage{natbib}
\newcommand{\eg}{e.\,g.\ }
\newcommand{\ie}{i.\,e.\ }

\begin{document}
\title{Generalized Fourier transform applied to fractional derivatives and Quantum Statistics Distributions.}
\author{Cyril Belardinelli}
\email[Electronic mail: ]{cyril.belardinelli@ksso.ch}
\affiliation{Kantonsschule Solothurn, Herrenweg 18, CH-4502 Solothurn, Switzerland}
\date{\today}
\begin{abstract}
In the present article the author extends the Fourier transform to a more general class of functions; First to power-law functions with integer and half-integer exponents then to the widely used quantum statistics function (Fermi-Dirac and Bose-Einstein-distributions). The results are used to derive rather straightforwardly an integral formula. Moreover, generalized Fourier transforms are used to define fractional derivatives of distributions (generalized functions) and Fourier series (periodic functions) as well.
\end{abstract}
\keywords{Fourier Transformation; temperate distributions; fractional derivative}
\maketitle
\section{\label{intro}Motivation}
In order to define the Fourier transform of a function in an elementary sense one has (at least) to require that it decays fast enough at infinity. This restriction is in many ways disturbing. In the past many attempts have been made to extend the definition to a more general class of functions.\cite{bochner:1932,carleman:1944,schwartz:1951, gelfand:1958} 
As is well-known the theory of distributions offers the possibility to extend the definition to functions in ${L}^p$-spaces (and beyond). 
The usual definition of a Fourier transform is given by:
\begin{equation}
\mathcal{F}(f)(k)=\hat{f}(k):=\int_{\mathbb{R}}f(x)e^{-ikx}\mathrm{d}x, \quad x\in \mathbb{R}
\end{equation}
provided the integral exists, \eg $|f(x)|$ is integrable in Lebesque's sense. Usually one restricts the definition to $L^{1}$-functions. However, the function $f(x)=\frac{1}{\sqrt{x}}$ does not belong to $L^{1}$ but the Fourier transform for non-vanishing $k$-values can be calculated (on a given branch of $\log{z}, z \in \mathbb{C}$). Moreover, the constant $\bold{1}$-function has a well-defined Fourier transform in the distributional sense: $\mathcal{F}[\bold{1}]=2\pi\delta(k)$. Examples of this kind suggest therefore a generalization of the concept. The present article is organized in the following manner: Sec.~\ref{Sec.1} is partly expository; Fourier transforms of power-law functions with integer and half-integer exponents are calculated in a refreshing manner by using Heaviside's unit step function. Sec.~\ref{Sec.2} shows that the Fourier transforms of the widely used quantum statistics distributions (\ie Fermi-Dirac and Bose-Einstein distributions) do exist in the distributional sense. This calculation appear to be novel. Sec.~\ref{sec_frac_deriv} shows that fractional derivatives of distributions (generalized functions) are well-defined and fullfill the usual properties of common derivatives. Some formulas such as the fractional derivative of Heaviside's step function are rederived. Moreover, fractional derivatives of Fourier series (\ie periodic functions) are  introduced.  
In Sec.~\ref{special_integral} the integrals $\mathbf{{\int_{-\infty}^{\infty}\frac{\sin^{n}x}{x^m}\mathrm{d}x}}$, $\mathbf{{n\geq m \geq 1} \in \mathbb{N}}$ are calculated in a remarkable straightforward manner by using the results of Sec.~\ref{Sec.1}. In App.~\ref{app} the more general case of $\mathbf{{\int_{0}^{\infty}\frac{\sin^{n}x}{x^m}\mathrm{d}x}}$ is calculated giving rise to an integral formula which is appearently unknown in the literature.       
\section{\label{Sec.1}Fourier Transforms of power functions}
\subsection{Fourier transforms of $\bold{f(x)=x^{n}, n\in \mathbb{N}_{0}}$}
In this section we calculate Fourier transforms of rapidly increasing functions such as power functions $f(x)=x^n$ ($n>0$) in a mainly expository but refreshing way. This extension is possible when $f(x)$ is interpreted as a distribution rather than a function. To achieve this purpose one starts with the Fourier transform of the Heaviside unit step function defined as:
\begin{equation}
\Theta(x)=\begin{cases}
1&\text{if $x>0$}\\
0&\text{if $x<0$}
\end{cases}
\end{equation}
Throughout this paper the use of $\Theta(x)$ is rather important. One starts by calculating the following Fourier transform: 
\begin{equation}
\label{calcul_1}
\mathcal{F}[x^{n}\Theta(x)](k)\equiv\int_{0}^{\infty}x^{n}e^{-ikx}\mathrm{d}x
\end{equation}
Which is evaluated by using the identity: 
\begin{equation}
\label{identity}
\mathcal{F}\left[x^{n}\phi\right](k)=i^{n}\frac{d^{k}}{dk^{n}}\hat{\phi}(k) 
\end{equation}
Where $\phi(x)$ is usually taken from the Schwartz-space $\mathcal{S}$ of fast (at a polynomial rate) decreasing functions. Here, one assumes that (\ref{identity}) is also valid for temperate distributions.
One gets then:
\begin{equation}
\eqref{calcul_1}=i^{n}\frac{d^{k}}{dk^{n}}\hat{\Theta}(k)
\end{equation}
The Fourier transform of $\Theta(x)$ is well-known and reads as follows:
\begin{equation}
\hat{\Theta}(k)=\pi \delta(k)+\frac{1}{i k}
\end{equation}
Therefore, one yields:
\begin{equation}
\mathcal{F}\left[x^{n}\Theta(x)\right](k)=i^{n}\pi \,\delta^{(n)}(k)+i^{n-1}n!\frac{(-1)^{n}}{k^{n+1}}
\end{equation}
where $\delta^{(n)}$ denotes the $n_{\text{th}}$ derivative of the $\delta$-function.\\
Using the simple identity $\Theta(x)+\Theta(-x)=1$ one can derive:
\begin{equation}
\label{FT_power_pos}
\begin{split}
&\mathcal{F}[x^{n}](k)=\mathcal{F}[x^{n}(\Theta(x)+\Theta(-x))](k)\\
&=\mathcal{F}[x^{n}(\Theta(x)](k)+(-1)^{n}\mathcal{F}[(-x)^{n}(\Theta(-x)](k)\\
&=2\pi i^{n}\delta^{(n)}(k)
\end{split}
\end{equation}
where $\delta^{(n)}(-k)=(-1)^{n}\delta^{(n)}(k)$ has been used.
\subsection{Fourier transform of $\bold{f(x)=x^{-n}, n\in \mathbb{N}}$}
Once again the unit step function $\Theta(x)$ shows its utility. One starts with its Fourier transform:
\begin{equation}
\hat{\Theta}(k)=\pi \delta(k)+\frac{1}{i k}
\end{equation}
considering the identity 
\begin{equation}
\Theta(x)=\frac{1}{2}+\frac{1}{2}\text{sgn}(x)
\end{equation}
one derives immediately:
\begin{equation}
\mathcal{F}[\text{sgn(x)}](k)=\frac{2}{ik}
\end{equation}
By Fourier transforming the last identity one concludes:
\begin{equation}
\mathcal{F}[x^{-1}](k)=-i\pi \text{sgn}(k)
\end{equation}
The latter result is nothing but an expression of Dirichlet's discontinous integral:
\begin{equation}
\int_{-\infty}^{\infty}\frac{e^{ikx}}{x}\mathrm{d}x=i\pi \text{sgn}(x)
\end{equation}
The Fourier transform of $f(x)=x^{-2}$ can be evaluated easily by using 
\begin{equation}
\frac{d}{\mathrm{d}x}x^{-1}=-x^{-2} \quad\text{and}\quad \mathcal{F}\left[\frac{d}{\mathrm{d}x}x^{-1}\right]=ik\mathcal{F}\left[x^{-1}\right]
\end{equation}
A short calculation yields
\begin{equation}
\mathcal{F}\left[x^{-2}\right]=k\pi\text{sgn}(k)
\end{equation}
In general one gets the formula:
\begin{equation}
\label{FT_power}
\mathcal{F}\left[x^{-n}\right](k)=\frac{\pi k^{n-1}}{i^{n}(n-1)!}\text{sgn}(k), \quad n \in \mathbb{N}
\end{equation}
\subsection{Fourier transforms of $\bold{f(x)=x^{n+\frac{1}{2}}, n\in \mathbb{Z}}$}
The purpose of this section is to calculate the Fourier transform of $f(x)=x^{n+\frac{1}{2}}$ for $n \in \mathbb{Z}$. One starts by calculating the integral:
\begin{equation}
\label{FT_power_half_integer}
\lim_{\epsilon\to\infty}\int_{0}^{\infty}x^{n}e^{-ikx-\epsilon x}=\lim_{\epsilon\to\infty}\frac{n!}{(\epsilon+ik)^{n+1}}=n! (ik)^{-n-1}
\end{equation} 
By using the Gamma-function $\Gamma(z)$ which is holomorphic in $z \in \mathbb{C}\setminus\{0,-1,-2,-3,...\}$ one can try to extend  (\ref{FT_power_half_integer}) to complex exponents:
\begin{equation}
\lim_{\epsilon\to\infty}\int_{0}^{\infty}x^{\alpha}e^{-ikx-\epsilon x}=\Gamma(\alpha+1) (ik)^{-\alpha-1}
\end{equation}
By restricting to half integer values $\alpha=n+\frac{1}{2}$ and extending the domain of integration to the entire real axis $\bold{\mathbb{R}}$ one \emph{defines} then:
\begin{equation}
\mathcal{F}[x^{n+\frac{1}{2}}](k):=\lim_{\epsilon\to\infty}\int_{-\infty}^{\infty}x^{n+\frac{1}{2}}e^{-ikx-\epsilon x}=2 \Gamma(n+\frac{3}{2}) (ik)^{-n-\frac{3}{2}}
\end{equation}
\section{Fourier Transforms of Quantum statistics functions}
\label{Sec.2}
As a further example of a non-trivial generalized Fourier transform we calculate the generalized Fourier transform of the Fermi-Dirac (Bose-Einstein) function $F_{\text{FD}}(x)=\frac{1}{e^{\beta x}+1}$ $\left(F_{\text{BE}}(x)=\frac{1}{e^{\beta x}-1}\right)$which is extensively used in various contexts of quantum statistics. The calculation presented here is, however, either new or, at least hard to find in the literature. The parameter $\beta=\frac{1}{T}$ is the inverse temperature $T$ (with $k_{B}$ the Boltzmann constant). For convenience the chemical potential $\mu$ is set to zero.
In the limit of vanishing temperature $T$ (\ie $\beta \rightarrow \infty$) one has:
\begin{equation}
\label{limes_beta}
\lim_{\beta \to \infty}F_{\text{FD}}(x)=\Theta(-x)
\end{equation}
The Fourier transform of $F_{\text{FD}}(x)$ is formally given by:
\begin{equation}
\label{FT_Fermi_Dirac}
\hat{F}_{\text{FD}}(k)=\int_{\mathbb{R}}\frac{e^{-ikx}}{1+e^{\beta x}}\mathrm{d}x
\end{equation}
The latter integral is divergent in the classical sense. However, when $F_{\text{FD}}(x)$ is interpreted as a distribution (generalized function) one can calculate its Fourier transform as shown in the following calculation.
By splitting the integral one gets:
\begin{equation}
\hat{F}_{\text{FD}}(k)=\int^{\infty}_{0}\frac{e^{ikx}}{1+e^{-\beta x}}\mathrm{d}x+\int^{\infty}_{0}\frac{e^{-ikx-\beta x}}{1+e^{-\beta x}}\mathrm{d}x
\end{equation}
By expanding both denominators in a geometric series:
\begin{equation}
\frac{1}{1+e^{-\beta x}}=\sum^{\infty}_{n=0}(-1)^{n}e^{-\beta nx}
\end{equation}
one yields:
\begin{equation}
\label{sum}
\hat{F}_{\text{FD}}(k)=\int^{\infty}_{0}e^{ikx}\mathrm{d}x+ \sum^{\infty}_{n=1}(-1)^{n}\int^{\infty}_{0}\left[e^{-x(\beta n-ik)}-e^{-x\left(\beta n +ik\right)}\right]\mathrm{d}x
\end{equation}
The first integral in (\ref{sum}) is divergent. It reflects the fact that the Fourier transform of $F_{\text{FD}}(x)=\frac{1}{1+e^{\beta x}}$ does not exist in a classical sense. The second term in (\ref{sum}) is well-defined. At this stage one has simply isolated the divergent part of the (classically divergent) Fourier integral (\ref{FT_Fermi_Dirac}). However, the integral $\int^{\infty}_{0}e^{ikx}\mathrm{d}x$ has a well-defined meaning as a distribution. It is the Fourier transform of $\Theta(-x)$. The second term in (\ref{sum}) can easily be calculated. One gets then:
\begin{equation}
\label{sum2}
\hat{F}_{\text{FD}}(k)=\hat{\Theta}(-k)+2ik\sum^{\infty}_{n=1}\frac{(-1)^n}{\beta^{2}n^{2}+k^2}
\end{equation}
The latter sum can be evaluated explicitly via the formula (See \eg Ref.\cite{gradshteyn:2014}):
\begin{equation}
\sum^{\infty}_{n=1}\frac{(-1)^n}{\alpha^{2}-n^2}=\left(\frac{\pi}{\sin{\alpha \pi}}-\frac{1}{\alpha}\right)\frac{1}{2\alpha},\quad \alpha \in \mathbb{C/Z}
\end{equation}
One has then:
\begin{equation}
\hat{F}_{\text{FD}}(k)=\hat{\Theta}(-k)+\frac{1}{ik}+\frac{i\pi}{\beta\sinh{\frac{\pi k}{\beta}}}
\end{equation}
By using $\hat{\Theta}(-k)=\pi \delta(k)-\frac{1}{ik}$ one gets:
\begin{equation}
\hat{F}_{\text{FD}}(k)=\pi \delta(k)+\frac{i\pi}{\beta\sinh{\frac{\pi k}{\beta}}}
\end{equation}
As expected, one yields $$\lim_{\beta \to \infty}\hat{F}_{\beta}(k)=\hat{\Theta}(-k)$$
A calculation along a similar track yields the Fourier transform of the Bose-Einstein-distribution:
\begin{equation} 
\hat{F}_{\text{BE}}(k)=\frac{\pi}{i\beta}\coth{\frac{\pi k}{\beta}}-\pi \delta(k)
\end{equation}
\section{Fractional derivatives}
\label{sec_frac_deriv}
The purpose of this section is to motivate the Fourier transform of distributions and of Fourier series. We begin with some general considerations.
As commonly known, temperate distributions $u:\mathcal{S}\rightarrow \mathbb{C}$ ($\mathcal{S}$: Schwartz-space) have well-defined derivatives of arbitrary order and can be Fourier transformed as well. They are therefore good candidates for an extension to derivatives of fractional order. One starts with the formula for the $n_{\text{th}}$ derivative
of the temperate distribution $u(\phi)$ given in the form:
\begin{equation}
\label{derivat_distri}
\partial^{n}{u}(\phi)=(-1)^{n}u(\partial^{n}\phi)
\end{equation}
where $\phi \in \mathcal{S}$\\
Formula (\ref{derivat_distri}) serves as a definition for derivatives of fractional order of a temperate distribution: 
\begin{equation}
\label{frac_derivat}
\partial^{\alpha}{u}(\phi)=(-1)^{|\alpha|}u(\partial^{\alpha}\phi), \quad \alpha \in \mathbb{R}_{+}
\end{equation}
In formula (\ref{frac_derivat}) one has to calculate the fractional derivative of Schwartz-function $\phi(x) \in \mathcal{S}$ which can be done via the identity:
\begin{equation}
\label{formula_frac_derivat}
\mathcal{F}(\partial^{\alpha}\phi)(k)=i^{|\alpha|}k^{\alpha}\mathcal{F}(\phi)(k) \in \mathcal{S}
\end{equation}
Backtransforming the right-hand side of (\ref{formula_frac_derivat}) yields $\partial^{\alpha}\phi(x)\, \in \mathcal{S}$. This means that
$u(\partial^{n}\phi)$ is again a temperate distribution. Since a very wide class of functions can be interpreted as temperate distributions which have been proved to be fractional derivable one has at least in principle a method at hand to calculate the fractional derivative of the function itself. \\
Since 
\begin{equation}
\mathcal{F}({\partial^{\alpha}\partial^{\beta}\phi(x)})=i^{|\alpha|}k^{\alpha}\mathcal{F}({\partial^{\beta}\phi(x)})=
i^{|\alpha|}i^{|\beta|}k^{\alpha}k^{\beta}\phi(k)=\mathcal{F}(\partial^{\alpha+\beta}\phi(x))
\end{equation}
the \emph{semi-group property} is preserved:
\begin{equation}
\partial^{\alpha}\partial^{\beta}\phi(x)=\partial^{\alpha+\beta}\phi(x)
\end{equation}
As a  \emph{$1^{\text{st}}$} example we consider the fractional derivative of $f(x)=e^{iax}$. Its Fourier transform reads:
\begin{equation}
\mathcal{F}\left[e^{iax}\right](k)=2\pi\delta(k-a)
\end{equation}
In order to get $\partial^{\alpha}e^{iax}$ one has to calculate the Fourier inverse.
\begin{equation}
\label{half_derivative}
\partial^{\alpha}e^{iax}=\mathcal{F}^{-1}\left(i^{|\alpha|}k^{\alpha}2\pi\,\delta(k-a)\right)=i^{|\alpha|}a^{\alpha}e^{iax}
\end{equation}
In particular for $\alpha=\frac{1}{2}$:
\begin{equation}
\partial^{\frac{1}{2}}e^{iax}=\sqrt{\frac{a}{2}}(1+i)(\cos{ax}+i\sin{ax})
\end{equation}
Presuming linearity of $\partial^{\alpha}$ one gets:
\begin{equation}
\begin{split}
\partial^{\frac{1}{2}}\cos{ax}&=\sqrt{\frac{a}{2}}(\cos{ax}-\sin{ax})\\
\partial^{\frac{1}{2}}\sin{ax}&=\sqrt{\frac{a}{2}}(\cos{ax}+\sin{ax})
\end{split}
\end{equation}
One can easily check that in the specific example above the \emph{semi-group property}  is fullfiled:
\begin{equation}
\label{half_deriv_sin_cos}
\begin{split}
&\partial^{\frac{1}{2}}\partial^{\frac{1}{2}}\cos{ax}=-\sin{ax}\\
&\partial^{\frac{1}{2}}\partial^{\frac{1}{2}}\sin{ax}=+\cos{ax}\\
\end{split}
\end{equation}
With formula (\ref{half_derivative}) one has a natural method at hand to calculate the fractional derivative of any (reasonable) periodic function. Next, we calculate fractional derivatives of elementary functions.
\subsubsection{fractional derivative of $\Theta(x)$ and $\delta(x)$}
Fractional derivatives of Heaviside's step function $\Theta(x)$ and Dirac's-delta function $\delta(x)$ are straightforwardly derivable. For the sake of illustration we restrict to fractional derivatives of order $\frac{1}{2}$. 
\begin{equation}
\begin{split}
&\partial^{\frac{1}{2}}\Theta(x)=\mathcal{F}^{-1}\left[\sqrt{ik}\hat{\Theta}(k)\right]\\
&=\mathcal{F}^{-1}\left[\sqrt{ik}(\pi \delta(k)+(ik)^{-1})\right]\\
&\equiv\mathcal{F}^{-1}\left[(ik)^{-\frac{1}{2}}\right]=\sqrt{\frac{1}{\pi x}}
\end{split}
\end{equation}
The latter result has been derived by O. Heaviside in a different manner.\cite{heaviside:1970}(See also Ref.\cite{courant:1943}). However, the formula was actually known long before Heaviside.\cite{ross:1977}
To calculate $\partial^{\frac{1}{2}}\delta(x)$ one can use $H^{\prime}=\delta$ and commutative property (of fractional derivatives):
\begin{equation}
\partial^{\frac{1}{2}}\delta(x)=\partial^{\frac{1}{2}}H^{\prime}(x)=\partial\partial^{\frac{1}{2}}\Theta(x)=\partial\sqrt{\frac{1}{\pi x}}
=-\frac{2}{\sqrt{\pi}}x^{-\frac{3}{2}}
\end{equation}
\subsubsection{fractional derivative of a constant function}
The fractional derivative of a constant function $f(x)=c$ with $c \in \mathbb{C}$ can be calculated in the following way:
\begin{equation}
\mathcal{F}[c](k)=2\pi c \,\delta(k)
\end{equation}
As expected, the Fourier inverse gives:
\begin{equation}
\partial^{\alpha}f(x)=\mathcal{F}^{-1}\left[2\pi i^{|\alpha|} k^{\alpha}\delta(k)\right]=0 \quad \forall \alpha>0
\end{equation}
where $k^{\alpha}\delta(k)=0$ for $\alpha>0$ has been used.
\subsubsection{power functions $f(x)=x^{n}$}
First, we treat the case of positive integer exponents $n \in \mathbb{N}$. In order to be a true extension, the fractional derivative's definition must obviously include the conventional one. We check this by calculating the $m_{\text{th}}$ derivative of $f(x)=x^{n}$.\\
According to formula \eqref{FT_power_pos} one has:
\begin{equation}
\mathcal{F}[x^{n}](k)=2\pi i^{n}\frac{d^{n}}{dk^{n}}\delta(k)=2\pi(-i)^{n}\frac{n!}{k^{n}}\delta(k)
\end{equation}
Therefore:
\begin{equation}
\begin{split}
\partial^{m}x^{n}&=\mathcal{F}^{-1}\left(i^{m}k^{m}2\pi(-i)^{n}\frac{n!}{k^{n}}\delta(k)\right)\\
&=\frac{n!}{m!}\mathcal{F}^{-1}\left(2\pi(-i)^{n-m}\frac{(n-m)!}{k^{n-m}}\delta(k)\right)\\
&=n(n-1)(n-2)\cdots(n-m+1)x^{n-m}
\end{split}
\end{equation}
\subsubsection{power functions $f(x)=x^{-n}$}
According to \eqref{FT_power}:
\begin{equation}
\begin{split}
&\mathcal{F}\left[\frac{d^{m}}{\mathrm{d}x^{m}}x^{-n}\right](k)=i^{m}k^{m}\mathcal{F}[x^{-n}](k)\\
&=\frac{\pi k^{n+m-1}}{i^{n-m}(n-1)!}
\end{split}
\end{equation}
Fourier inverting the latter expression gives:
\begin{equation}
\frac{d^{m}}{\mathrm{d}x^{m}}x^{-n}=(-1)^{m}\frac{(n+m-1)!}{(n-1)!}x^{-(n+m)}
\end{equation}
\subsection{fractional derivative of Fourier series}
As already mentioned in Sec.~\ref{sec_frac_deriv} one can define through Eq.~(\ref{half_derivative}) fractional derivatives of Fourier series in a natural way.
In the following we give two examples.\\
As a first example we calculate the fractional derivative of order $\frac{1}{2}$ of the $2\pi$-periodical function:
\begin{equation}
f(x)=x \quad \text{for}-\pi<x<\pi \quad \text{(periodically continued)}
\end{equation}
Its Fourier series is given by:
\begin{equation}
f(x)=2\sum_{n=1}^{\infty}(-1)^{n-1}\frac{\sin nx}{n}
\end{equation}
By applying Eq.~(\ref{half_deriv_sin_cos}) one gets:
\begin{equation}
\label{weierstrass}
\partial^{\frac{1}{2}}f(x)=\sqrt{2}\sum_{n=1}^{\infty}(-1)^{n-1}\frac{\sin{nx}+\cos{nx}}{\sqrt{n}}
\end{equation}
\begin{figure}
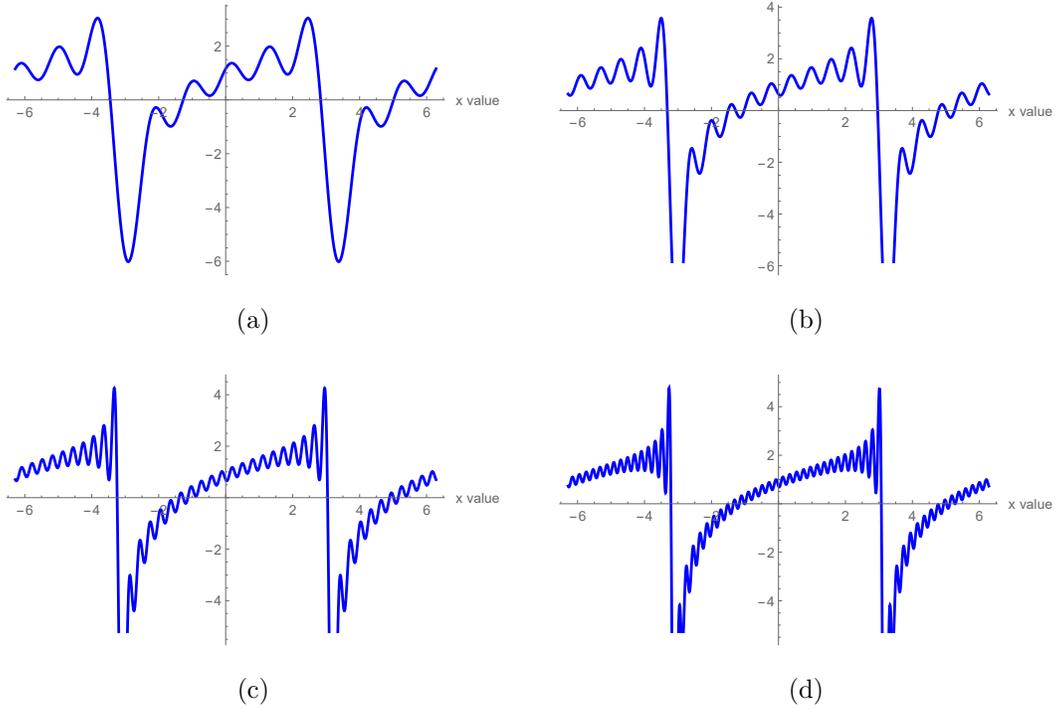
%
  \centering
  \subfloat[][]{\includegraphics[width=0.4\linewidth]{Bild_half_derivat_5.png}}%
  \qquad
  \subfloat[][]{\includegraphics[width=0.4\linewidth]{Bild_half_derivat_10.png}}%
 \qquad
  \subfloat[][]{\includegraphics[width=0.4\linewidth]{Bild_half_derivat_20.png}}%
 \qquad
  \subfloat[][]{\includegraphics[width=0.4\linewidth]{Bild_half_derivat_30.png}}%
  \caption{Fourier series of $\partial^{\frac{1}{2}}f(x)$ up to order: $n=5$ (a), $n=10$ (b), $n=20$ (c), $n=30$ (d).}%
\label{Bild}
\end{figure}
As shown in  Fig.~\ref{Bild} the truncated Fourier series is oscillating more rapidly with increasing order $n$. One can therefore conjecture that the Fourier series (\ref{weierstrass}) represents a function of the Weierstrass type with the property of being continous but nowhere differentiable.\cite{weierstrass:1895}
The second example is given by the function:
\begin{equation}
g(x)=|x| ,\quad\text{for}-\pi<x<\pi \quad \text{(periodically continued)}
\end{equation}
which has the Fourier series:
\begin{equation}
g(x)=\frac{\pi}{2}-\frac{4}{\pi}\sum_{n=1}^{\infty}\frac{\cos{(2n+1)x}}{(2n+1)^2}
\end{equation}
Once again, by Eq.~(\ref{half_deriv_sin_cos}) one has:
\begin{equation}
\label{2nd_example}
\partial^{\frac{1}{2}}g(x)=-\frac{4}{\pi\sqrt{2}}\sum_{n=1}^{\infty}\frac{\cos{(2n+1)x}-\sin{(2n+1)x}}{\sqrt{2n+1}(2n+1)}
\end{equation}
\begin{figure}[h!]
\centering
\includegraphics[width=8cm]{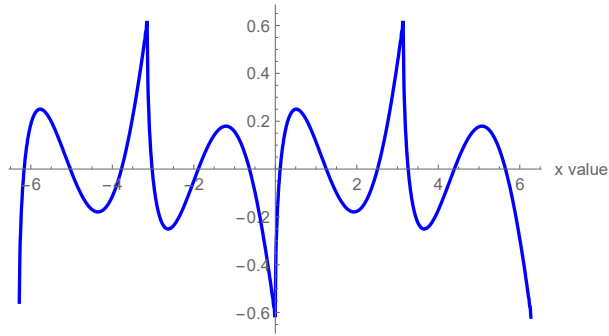}
\caption{Fourier series of $\partial^{\frac{1}{2}}g(x)$ up to order: $n=100$}
\end{figure}
\newpage
\section{Application to an Integral}
\label{special_integral}
\subsection{The integrals $\bold{\int_{\mathbb{R}}\frac{\sin^{n}x}{x^m}\mathrm{d}x}$ and $\bold{\int_{0}^{\infty}\frac{\sin^{n}x}{x^m}\mathrm{d}x}$   for $\bold{n\geq m \geq 1} \in \mathbb{N}$}
Formula (\ref{FT_power}) of the last section offers a nice and very short way of calculating integrals of the following type ($n\geq m$):
\begin{equation} 
\int_{\mathbb{R}} \frac{\sin^{n}x}{x^m}\mathrm{d}x, \quad n,m \in \mathbb{N}
\end{equation}
Historically important cases are the Dirichlet-Integral $\int^{\infty}_{0}\frac{\sin x}{x}\mathrm{d}x=\frac{\pi}{2}$ 
and $\int^{\infty}_{0}\frac{\sin^{2} x}{x^2}\mathrm{d}x=\frac{\pi}{2}$. The latter integral is used to prove the theorem of Wiener and Ikehara which is the basis for the shortest proof of the prime number theorem. One starts by writing:
\begin{equation}
\label{class_integrals}
\int_{\mathbb{R}} \frac{\sin^{n}x}{x^m}\mathrm{d}x=\frac{1}{2^{n}i^{n}}\int_{\mathbb{R}}\frac{\left(e^{ix}-e^{ix}\right)^{n}}{x^{m}}\mathrm{d}x
\end{equation}
By expanding the brackets and reversing summation and integration one gets:
\begin{equation}
(\ref{class_integrals})=\frac{1}{2^{n}i^{n}}\sum^{n}_{l=0}(-1)^{l}\binom{n}{l}\int_{\mathbb{R}}\frac{e^{-ix(2l-n)}}{x^m}\mathrm{d}x 
\end{equation}
The latter integral is classically divergent. However, by using formula (\ref{FT_power}) one gets immediately:
\begin{equation}
\label{formula_integral}
\int_{\mathbb{R}} \frac{\sin^{n}x}{x^m}\mathrm{d}x=\frac{\pi}{2^{n}i^{n+m}(m-1)!}\sum^{n}_{l=0}(-1)^{l}\binom{n}{l}(2l-n)^{m-1}\text{sgn}(2l-n) 
\end{equation}
Evidently, one gets zero if one of the exponents $n,m$ is odd.
A few examples:
\begin{equation}
\begin{split}
\int_{\mathbb{R}} \frac{\sin^{3}x}{x^3}\mathrm{d}x&=\frac{3}{4}\pi \quad \int_{\mathbb{R}} \frac{\sin^{4}x}{x^4}\mathrm{d}x=\frac{2}{3}\pi\\
\int_{\mathbb{R}} \frac{\sin^{5}x}{x^5}\mathrm{d}x&=\frac{115}{192}\pi \quad \int_{\mathbb{R}} \frac{\sin^{6}x}{x^6}\mathrm{d}x=\frac{11}{20}\pi
\end{split}
\end{equation}
In the case $n=m$ one gets by a simple rearrangement of the sum (\ref{formula_integral}) an equivalent formula:
\begin{equation}
\int_{\mathbb{R}} \frac{\sin^{n}x}{x^n}\mathrm{d}x=\frac{\pi}{2^{n-1}(n-1)!}\sum^{n/2}_{l=0}(-1)^{l}\binom{n}{l}(n-2l)^{n-1} 
\end{equation}
The latter formula is elegantly derived in Ref.\cite{apostol:1980} in another way.
It is possible to derive a formula for the more general case (See Appendix \ref{app}):
\begin{equation}
\label{formula}
\int_{0}^{\infty} \frac{\sin^{n}x}{x^m}\mathrm{d}x=\frac{i^{m-n+1}}{2^{n}(m-1)!} \sideset{}{'}\sum^{n}_{l=0}(-1)^{l}\binom{n}{l}(n-2l)^{m-1}\ln[{i(2l-n)}]
\end{equation}
where $\sideset{}{'}\sum$ means that the term $2l=n$ is omitted in (\ref{formula}). (Alternatively one could also use the convention $"0 \ln 0=0"$ which is meaningful since $\lim_{x\to 0} x\ln{x}=0$).\\
Once again a few examples:
\begin{equation}
\begin{split}
\int_{0}^{\infty} \frac{\sin^{3}x}{x^2}\mathrm{d}x&=\frac{3}{4}\ln{3} \quad \int_{0}^{\infty} \frac{\sin^{4}x}{x^3}\mathrm{d}x=\ln{2}\quad \int_{0}^{\infty} \frac{\sin^{3}x}{x}\mathrm{d}x=\frac{\pi}{4}\pi \\
\int_{0}^{\infty} \frac{\sin^{5}x}{x^4}\mathrm{d}x&=\frac{125}{96}\ln{5}-\frac{45}{32}\ln{3}
\end{split}
\end{equation}
\appendix
\section{The integrals $\bold{\int_{0}^{\infty}\frac{\sin^{n}x}{x^m}\mathrm{d}x}$, $\bold{n\geq m \geq 1} \in \mathbb{N}$}
\label{app}
To derive formula (\ref{formula}) one starts with the the integral:
\begin{equation}
f(k)=\int_{0}^{\infty}\frac{\sin^{n}x}{x^m}e^{-kx}\mathrm{d}x
\end{equation}
In the following we will use the fact that $\lim_{k\to \infty} f(k)=0$.\\
The $m_{th}$ derivative reads as follows:
\begin{equation}
\label{formula_general}
\begin{split}
&f^{(m)}(k)=(-1)^{m}\int_{0}^{\infty}\sin{^{n}x}e^{-kx}\mathrm{d}x\\
&=\frac{(-1)^{m}}{(2i)^{n}}\sum^{n}_{l=0}(-1)^{l}\binom{n}{l}\int_{0}^{\infty}  e^{-ix(2l-n)-kx}\, \mathrm{d}x\\
&=\frac{(-1)^{m}}{(2i)^{n}}\sum^{n}_{l=0}\binom{n}{l}\frac{(-1)^l}{k+i(2l-n)}
\end{split}
\end{equation}
In order to regain the function $f(k)$ one has to calculate the $m_{th}$ antiderivative of the latter expression. To do so we use the simple fact that the $m_{th}$ antiderivative of $g(y)=\frac{1}{y}$ has the form:
\begin{equation}
\label{expression}
\int\mathrm{d}y\int\mathrm{d}y\cdots\int\mathrm{d}y\,\frac{1}{y}=A_{m-1}y^{m-1}\log{y}-B_{m-1}y^{m-1}+C_{1}+C_{2}y+\dots+C_{m}y^{m-1} 
\end{equation}
Where $A_{m-1}$ and $B_{m-1}$ are some coefficient depending on $m$. $C_1, C_2 \dots C_m$ are integration constants. One derives easily a recursion formula for the coefficients $A_{m-1}$ and $B_{m-1}(A_0=B_0=1)$:
\begin{equation}
\begin{split}
A_{m-1}&=\frac{A_{m-2}}{m-1}\\
B_{m-1}&=\frac{B_{m-2}}{m-1}+\frac{A_{m-2}}{(m-1)^2}
\end{split}
\end{equation}
From which one gets $A_{m-1}=\frac{1}{(m-1)!}$. As explained below, an explicit formula for $B_{m-1}$ is not needed. 
Next, one sets $y=k+i(2l-n)$ and inserts (\ref{expression}) into (\ref{formula_general}). \\Since one has evidently
\begin{equation}
\lim_{k\to \infty} f^{(l)}(k)=0, \quad \text{for}\quad l=0,1\dots,m
\end{equation}
the polynomial terms appearing in the sum (\ref{formula_general}) cancel each other. It is therefore enough to consider only the coefficients  $A_{m-1}=\frac{1}{(m-1)!}$. Doing so one yields the formula:
\begin{equation}
f(0)=
\lim_{k\to 0} f(k)=
\frac{i^{m-n+1}}{2^{n}(m-1)!} \sideset{}{'}\sum^{n}_{l=0}(-1)^{l}\binom{n}{l}(n-2l)^{m-1}\log[{i(2l-n)}]
\end{equation}
By using the identity $\log[{i(2l-n)}]=\log{|2l-n|}+i\frac{\pi}{2}\text{sgn}(2l-n)$ one can write:
\begin{equation}
\begin{split}
\int_{0}^{\infty}\frac{\sin^{n}x}{x^m}\mathrm{d}x &=\frac{i^{m-n+1}}{2^{n}(m-1)!} \sideset{}{'}\sum^{n}_{l=0}(-1)^{l}\binom{n}{l}(n-2l)^{m-1}\log{|2l-n|}\\
&-\frac{i^{m-n}}{2^{n}(m-1)!}\frac{\pi}{2} \sideset{}{'}\sum^{n}_{l=0}(-1)^{l}\binom{n}{l}(n-2l)^{m-1}\text{sgn}(2l-n)
\end{split}
\end{equation}
The latter sum vanish when $m-n$ is odd. Whereas the former sum vanish when $m-n$ is even. Therefore, one gets the general formulae :
\begin{equation}
\int_{0}^{\infty}\frac{\sin^{n}x}{x^m}\mathrm{d}x = \left\{
\begin{array}{lll}
-\frac{i^{m-n}}{2^{n}(m-1)!}\frac{\pi}{2} \sideset{}{'}\sum^{n}_{l=0}(-1)^{l}\binom{n}{l}(n-2l)^{m-1}\text{sgn}(2l-n) &\textrm{if $n-m$ even} \\
& \\
\frac{i^{m-n+1}}{2^{n}(m-1)!} \sideset{}{'}\sum^{n}_{l=0}(-1)^{l}\binom{n}{l}(n-2l)^{m-1}\log{|2l-n|} & \, \textrm{if $n-m$ odd} \\
\end{array}
\right. 
\end{equation}
Obviously, in both cases one gets real values (absence of imaginary terms).
\appendix*
\begin{acknowledgments}
The author thanks Thomas Petermann for inspiring discussions.
\end{acknowledgments} 
\appendix*
\newpage
\bibliography{References_generalized_fourier}
\bibliographystyle{plain}
\end{document}